\documentclass[twoside, final]{amsart} 

\usepackage{amsmath, amsfonts, amsthm, amssymb}
\usepackage{exscale}
\usepackage{cite}
\usepackage{epsfig}
\usepackage{amscd}
\usepackage{graphics}
\usepackage{color}
\usepackage{dsfont}

\usepackage[notref, notcite]{showkeys}

\renewcommand{\geq}{\geqslant}
\renewcommand{\leq}{\leqslant}

\newtheorem{theorem}{Theorem}

\newtheorem{proposition}{Proposition}[section]

\newtheorem{corollary}[proposition]{Corollary}
\newtheorem{lemma}[proposition]{Lemma}
\newtheorem*{main-theorem}{Main Theorem}
\newtheorem*{theorem*}{Theorem}

\theoremstyle{definition}

\newtheorem{remark}[proposition]{Remark}
\newtheorem{definition}[proposition]{Definition}
\newtheorem*{remark*}{Remark}

\numberwithin{equation}{section}

\def\phi{\varphi}

\def\ZZ{{\mathbb Z}}

\def\reals{{\mathbb R}}

\def\Ci{{\mathcal C}^\infty}

\def\WF{\mathrm{WF}_h}

\def\supp{\mathrm{supp}\,}

\def\O{{\mathcal O}}
\def\SS{{\mathbb S}}
\def\s{{\mathcal S}}

\def\phi{\varphi}
\def\half{{\frac{1}{2}}}

\def\be{\begin{eqnarray*}}
\def\ee{\end{eqnarray*}}
\def\ben{\begin{eqnarray}}
\def\een{\end{eqnarray}}

\def\lll{\left\langle}
\def\rrr{\right\rangle}

\def\L2R{L_{\text{Rest}}^2}

\def\11{\mathds{1}}

\def\tchi{\tilde{\chi}}
\def\L2c{L^2_{\text{comp}}}

\def\tDelta{\widetilde{\Delta}}

\def\tu{\tilde{u}}

\def\tA{\tilde{A}}

\def\Vol{\text{Vol}}

\def\tE{\widetilde{E}}
\def\p{\partial}
\def\GG{\mathcal{G}}

\begin{document}

\title[Unique Continuation]{Unique Continuation for
   Quasimodes 
   on Surfaces of Revolution: Rotationally invariant Neighbourhoods}
\author{Hans Christianson}
\address{Department of Mathematics, University of North Carolina}
\email{hans@math.unc.edu}


 \thanks{The author is extremely grateful to Andr\'as Vasy and Jared
   Wunsch for their help with the moment map definition of irreducible
   quasimode, and for their comments and suggestions on earlier versions of this
   paper.  
This work is supported in part by NSF grant 
   DMS-0900524.}

\subjclass[2000]{}

\begin{abstract}
We prove a strong conditional unique continuation estimate for
irreducible quasimodes in rotationally invariant neighbourhoods on compact surfaces of revolution.  The
estimate states that Laplace quasimodes which cannot be decomposed
as a sum of other quasimodes have $L^2$ mass bounded below by
$C_\epsilon \lambda^{-1 - \epsilon}$ for any $\epsilon>0$  on any
open rotationally invariant neighbourhood which meets the semiclassical wavefront set of the quasimode.
For an analytic manifold, we conclude the same estimate with a lower
bound of $C_\delta \lambda^{-1 + \delta}$ for some fixed $\delta>0$.

\end{abstract}

\maketitle

\section{Introduction}
\label{S:intro}

We consider a compact periodic surface of revolution $X = \SS^1_x \times \SS^1_\theta$, equipped with a metric of the form
\[
ds^2 = d x^2 + A^2(x) d \theta^2,
\]
where $A \in \Ci$ is a smooth function, $A \geq \epsilon>0.$  Our
analysis is microlocal, so applies also to any compact surface of revolution with no
boundary, and to certain surfaces of revolution with boundary under
mild assumptions, however we will concentrate on the toral case for
ease of exposition.



From such a metric, we get the volume form
\[
d \Vol = A(x) dx d \theta,
\]
and the Laplace-Beltrami operator acting on $0$-forms
\[
\Delta f = (\partial_x^2 + A^{-2} \partial_\theta^2 + A^{-1}
A' \partial_x) f.
\]

We are concerned with quasimodes, which are the building blocks from
which eigenfunctions are made, however we need to define the most
basic kind of quasimodes, which we will call {\it irreducible
  quasimodes}, meaning the quasimodes which cannot be decomposed as a
sum of two or more nontrivial quasimodes.  In order to make our
definitions, we recall first that the geodesic flow on $T^*X$
is the Hamiltonian system associated to the principal symbol of the
Laplace-Beltrami operator:
\[
p(x, \xi, \theta, \eta) = \xi^2 + A^{-2}(x) \eta^2.
\]
A fixed energy level $p = \text{ const.}$ consists of all the
geodesics of that constant ``speed''.  
For the case of the geodesic Hamiltonian system on
$T^*X$, there are two conserved quantities, the total energy and the
angular momentum $\eta^2$.  The {\it moment map} is the map sending
points of $T^*X$ to their associated conserved quantities, that is
\[
M(x, \xi, \theta, \eta) = \left( \begin{array}{c} \xi^2 + A^{-2}(x)
    \eta^2 \\ \eta^2 \end{array} \right).
\]
When the gradient of $M$ has rank $2$, then $M$ defines a submersion,
so each connected component of the preimage is a $2$-manifold.
Points in $T^*X$ where $M$ has rank $1$ or $0$ are called critical points,
and points
$(P,Q) \in \reals^2$ such that $\{ M = (P,Q) \}$ contains critical
points are called critical values.  Critical points 
 correspond to
latitudinal periodic geodesics, which can also carry quasimode mass,
and critical values have preimages which may have infinitely many
latitudinal periodic geodesics.  
The semiclassical wavefront set is always a closed invariant subset of the
energy surface, so our definition of irreducible quasimode will be one which has wavefront mass confined to the closure of one of
these two kinds of sets, distinguished by rank of $M$.

\begin{definition}
An {\it irreducible} quasimode is a quasimode whose semiclassical
wavefront set is contained in the closure of a single connected
component in $T^*X$ where the moment map has constant rank.

\end{definition}

\begin{figure}
\hfill
\centerline{\input{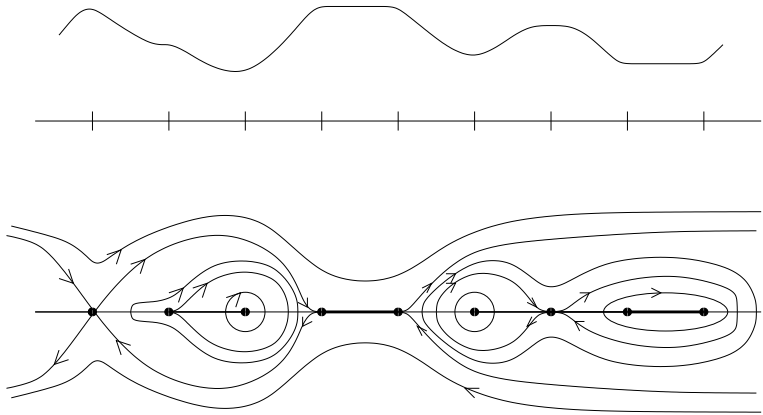}}
\caption{\label{fig:ex-pot} The reduced phase space of a toral surface of revolution with many
  periodic latitudinal geodesics.}
\hfill
\end{figure}

We also will require a limit on the geodesic complexity by assuming
there are only a finite number of connected regions of latitudinal
periodic geodesics.  This will not preclude having infinitely many
periodic latitudinal geodesics, but merely having accumulation points
of connected components of latitudinal geodesics.  We therefore will
assume that the moment map has a finite number of critical values,
each of which has a preimage of finitely many non-empty connected
components.  Note this allows intervals of latitudinal periodic
geodesics, but does not allow accumulation of such sets.  For an
example, see Figure \ref{fig:ex-pot}.

Finally, we will require a certain $0$-Gevrey regularity on the manifold,
which in a sense says our manifold is not too far from being
analytic.  
Such a $0$-Gevrey assumption nevertheless allows for
non-trivial functions which are constant on intervals, so this is a
very general class of manifolds.  Of course this includes analytic
manifolds, for which we have a stronger estimate.   See Subsection \ref{SS:gevrey} for the precise definitions.

\begin{theorem}
\label{T:T1}

Let $X$ be as above, for a generating curve in the $0$-Gevrey class $A(x)
\in \GG^0_\tau( \reals)$ for some $\tau < \infty$.  Assume the moment
map has finitely many critical values, with preimages consisting of
finitely many connected components.  Suppose $u$ is a (weak)
irreducible quasimode satisfying $\| u \| = 1$ and 
\[
(-\Delta - \lambda^2) u = \O( \lambda^{-\beta_0} ),
\]
for some fixed $\beta_0 >0$.  Let $\Omega \subset
X$ be a rotationally invariant neighbournood, $\Omega = (a,b)_x \times \SS_\theta^1$.  Then either 
\begin{enumerate}
\item
\[
\| u \|_{L^2(\Omega)} = \O( \lambda^{-\infty} ), 
\]
or 
\item
for any $\epsilon>0$, there exists $C = C_{\epsilon, \Omega, \beta_0}>0$ such that
\begin{equation}
\label{E:lower-bound}
\| u \|_{L^2(\Omega)} \geq C \lambda^{-1 - \epsilon}.
\end{equation}
\end{enumerate}

\end{theorem}

\begin{remark}

The proof will show that a more or less straightforward
commutator/contradiction argument gives a lower bound of
$\lambda^{-1-\beta_0}$.  The difficulty comes in trying to beat this
lower bound.

\end{remark}

In the analytic category, we have a significant improvement.  Of
course in the case of an analytic manifold, there can be no infinitely
degenerate critical elements, nor can there be any accumulation points
of sets of 
latitudinal periodic geodesics, so we do not need to make the
assumption about finite geodesic complexity.

\begin{corollary}
\label{C:C1}
Let $X$ be as above, and assume $X$ is analytic.   Suppose $u$ is a (weak)
irreducible quasimode satisfying $\| u \| = 1$ and 
\[
(-\Delta - \lambda^2) u = \O( 1 ).
\]
Then for any open rotationally invariant neighbourhood $\Omega \subset
X$, either 
\begin{enumerate}
\item
\[
\| u \|_{L^2(\Omega)} = \O( \lambda^{-\infty} ), 
\]
or 
\item
there exists a fixed $\delta>0$ and a constant $C = C_{\Omega}>0$ such that
\[
\| u \|_{L^2(\Omega)} \geq C \lambda^{-1 + \delta}.
\]
\end{enumerate}

\end{corollary}

\begin{figure}
\hfill
\centerline{\input{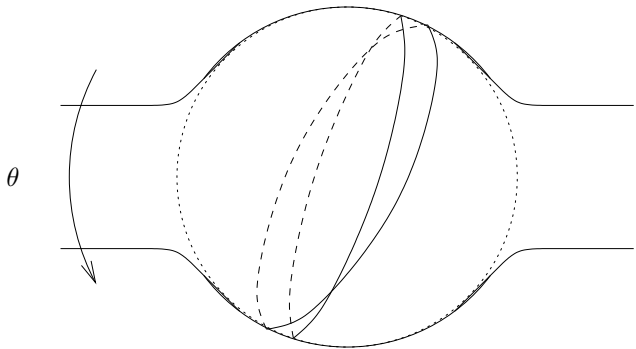}}
\caption{\label{fig:sph} A surface of revolution with a piece of
  $\SS^2$ embedded.  Also sketched are two ``isoenergetic'' periodic
  geodesics which are $\theta$ rotations of each other.  One can
  construct pathological quasimodes which are continuous, compactly
  supported 
  superpositions of isoenergetic quasimodes associated to such geodesics.}
\hfill
\end{figure}

\begin{remark}
The assumption that $\Omega \subset X$ is a rotationally invariant neighbourhood of
the form $\Omega = (a,b)_x \times \SS_\theta^1$ is necessary for this level
of generality.  To see this, consider the case where $X$ has part of a
2-sphere embedded in it.  Then there are many periodic geodesics close
to the latitudinal one.  But these geodesics can be rotated in
$\theta$ without changing the angular momentum.  Each one of these is
elliptic and can carry a Gaussian beam type quasimode.  Hence one can
create an irreducible quasimode as a superposition of these Gaussian
beams.  The resulting ``band'' of quasimodes need not have nontrivial
mass except in a rotationally invariant neighbourhood.  See Figure \ref{fig:sph}.

\end{remark}

\section{Preliminaries}
In this section we review some of the definitions and preliminary
computations necessary for Theorem \ref{T:T1}, as well as recall the
spectral estimates we will be using.

\subsection{The $0$-Gevrey class of functions}
\label{SS:gevrey}
For this paper, we use the
following 
$0$-Gevrey classes of functions with respect to order of vanishing,
introduced in \cite{Chr-inf-deg}.

\begin{definition}

For $0 \leq \tau <
\infty$, let $\GG^0_\tau (
\reals )$ be the set of all smooth functions $f : \reals \to \reals$
such that, for each $x_0 \in \reals$, there exists a neighbourhood $U \ni
x_0$ and a constant $C$ such that, for all $0 \leq s \leq k$, 
\[
| \partial_x^k f(x) -\partial_x^k f(x_0) | \leq C (k!)^C | x - x_0
|^{-\tau (k-s) } | \partial_x^s f(x) - \partial_x^s f(x_0) |, \,\,\, x
\to x_0 \text{ in } U.
\]

\end{definition}

This definition says that the order of vanishing of derivatives of a
function is only polynomially worse than that of lower derivatives.
Every analytic function is in one of the 0-Gevrey classes $\GG^0_\tau$ for some $\tau < \infty$, but many more functions are as well. 
For example, the function
\[
f(x) = \begin{cases} \exp (-1/x^p), \text{ for } x >0, \\ 0, \text{ for }x \leq 0
\end{cases}
\]
is in $\GG^0_{p+1}$, but
\[
f(x) = \begin{cases} \exp (-\exp(1/x)), \text{ for } x >0, \\ 0, \text{ for }x \leq 0
\end{cases}
\]
is not in any $0$-Gevrey class for finite $\tau$.

\subsection{Conjugation to a flat problem}

We observe that we can conjugate
$\Delta$ by an isometry of metric spaces and separate variables so
that spectral analysis of $\Delta$ is equivalent to a one-variable
semiclassical problem with potential.  That is, let $T : L^2(X, d
\Vol) \to L^2(X, dx d \theta)$ be the isometry given by
\[
Tu(x, \theta) = A^{1/2}(x) u(x, \theta).
\]
Then $\tDelta = T \Delta T^{-1}$ is essentially self-adjoint on $L^2 (
X, dx d \theta)$.  A simple calculation gives
\[
-\tDelta f = (- \partial_x^2 - A^{-2}(x) \partial_\theta^2 + V_1(x) )
f,
\]
where the potential
\[
V_1(x) = \frac{1}{2} A'' A^{-1} - \frac{1}{4} (A')^2 A^{-2}.
\]

If we now separate variables and write $\psi(x, \theta) = \sum_k
\phi_k(x) e^{ik \theta}$, we see that
\[
(-\tDelta- \lambda^2) \psi = \sum_k e^{ik \theta} P_k \phi_k(x),
\]
where
\[
P_k \phi_k(x) = \left(-\frac{d^2}{dx^2} + k^2 A^{-2}(x) + V_1(x) -
  \lambda^2 \right)
\phi_k(x).
\]
Setting $h = |k|^{-1}$ and rescaling, we have the semiclassical operator
\begin{equation}
\label{E:separation}
P(z,h) \phi(x) = (-h^2 \frac{d^2}{dx^2} + V(x) -z) \phi(x),
\end{equation}
where the potential is
\[
V(x) = A^{-2}(x) + h^2 V_1(x)
\]
and the spectral parameter is $z = h^2 \lambda^2$.  In Section
\ref{S:proofs} we will at first let $h = \lambda^{-1}$ be our
semiclassical parameter for the whole quasimode, but then switch to $h
= |k|^{-1}$ to estimate the parts of the quasimode microsupported
where the critical elements are located.  The relevant microlocal
estimates near critical elements are summarized in the following Subsection.

\subsection{Spectral estimates for weakly unstable critical sets}

In this subsection we summarize the spectral estimates we will use for weakly
unstable critical elements obtained in \cite{Chr-NC,Chr-NC-erratum,Chr-QMNC,ChWu-lsm,ChMe-lsm,Chr-inf-deg}.

\begin{definition}
Let $(P,Q)$ be a critical value of the moment map.  Then there are
points in $M^{-1}(P,Q)$ where the moment map has rank $1$ (or $0$, but
these points are easy to handle (see below)).  For these points, there are 
latitudinal periodic geodesics.  If the principal part of the
potential, $A^{-2}(x)$, for the reduced Hamiltonian $\xi^2 +
A^{-2}(x)$ has an ``honest'' minimum at $x_0$ in the sense that if
$[a,b]$ is the maximal closed interval containing $x_0$ with 
$
A^{-2}(x) = A^{-2}(x_0)$ on it, then $(A^{-2})' <0$ for $x < a$ in
some small neighbourhood, and $(A^{-2} )' >0$ for $x >b$ in some other
small neighbourhood, then we say this critical element is {\it weakly
  stable}.  In all other cases, we say the critical element is {\it
  weakly unstable}.

\end{definition}

In the following subsections, we review the microlocal estimates from
\cite{Chr-inf-deg} for weakly unstable critical elements.  Taken
together, they imply the following theorem.
\begin{theorem}
\label{T:inf-deg-est}
Let $\Lambda$ be a weakly unstable critical element in the reduced
phase space $T^*\SS^1_x$, and
assume $u$ has $h$-wavefront set sufficiently close
$\Lambda$.  Then for any $\epsilon>0$, there exists $C = C_\epsilon$
such that
\[
\| u \| \leq C h^{-2-\epsilon} \| ((hD)^2 + V(x) -z) u \|,
\]
for any $z \in \reals$.
\end{theorem}

\subsubsection{Unstable nondegenerate critical elements}

\label{SS:unst-nd}

A nondegenerate unstable critical element exists where the principal
part of the potential $V_0(x) =
A^{-2}(x)$ has a nondegenerate maximum.  To say that $x = 0$ is a nondegenerate maximum means that
 $x = 0$ is a critical point of
$V_0(x)$ satisfying $V_0'(0) = 0$, $V_0''(0) < 0$.

The following result as stated can be read off from \cite{Chr-NC,Chr-QMNC}, and has also been
studied in slightly different contexts in \cite{CdVP-I,CdVP-II} and
\cite{BuZw-bb}, amongst many others.  

\begin{lemma}
\label{L:ml-inv-0}
Suppose $x = 0$ is a nondegenerate local maximum of the principal part
of the potential $V_0$,
$V_0(0) = 1$.  For 
$\epsilon>0$ sufficiently small, let $\phi \in \s(T^* \reals)$
have compact support in $\{ |(x,\xi) |\leq \epsilon\}$.  Then there
exists $C_\epsilon>0$ such that 
\begin{equation}
\label{E:ml-inv-0}
\| P(z,h) \phi^w u \| \geq C_\epsilon \frac{h}{\log(1/h)} \|
\phi^w u \|, \,\,\, z \in [1-\epsilon, 1 + \epsilon].
\end{equation}
\end{lemma}

\begin{remark}
This estimate is known to be sharp, in the sense that the logarithmic
loss cannot be improved (see, for example, \cite{CdVP-I}).
\end{remark}

\subsubsection{Unstable finitely degenerate critical elements}

\label{SS:unst-fin}

In this subsection, we consider an isolated critical point at an 
unstable but finitely degenerate maximum.  That is, we now assume
that $x = 0$ is a degenerate maximum for the function $V_0(x) = A^{-2}(x)$ of order
$m \geq 2$.  If we again assume $V_0(0) = 1$, then this means that near
$x = 0$, $V_0(x) \sim 1 - x^{2m}$.  Critical points of this form were
first 
studied in \cite{ChWu-lsm}.

This Lemma and the proof are given in \cite[Lemma 2.3]{ChWu-lsm}.  

\begin{lemma}
\label{L:ml-inv-1}
For $\epsilon>0$ sufficiently small, let $\phi \in \s(T^* \reals)$
have compact support in $\{ |(x,\xi) |\leq \epsilon\}$.  Then there
exists $C_\epsilon>0$ such that 
\begin{equation}
\label{E:ml-inv-1}
\| P(z,h) \phi^w u \| \geq C_\epsilon h^{2m/(m+1)} \|
\phi^w u \|, \,\,\, z \in [1-\epsilon, 1 + \epsilon].
\end{equation}
\end{lemma}

\begin{remark}
This estimate is known to be sharp, in the sense that the exponent
$2m/(m+1)$ cannot be improved (see \cite{ChWu-lsm}).
\end{remark}

\subsubsection{Finitely degenerate inflection transmission critical elements}

We next study the case when the principal part of the potential has an inflection point of
finitely degenerate type.  That is, let us assume the point $x = 1$ is
a finitely degenerate inflection point, so that locally near $x = 1$,
the potential $V_0(x) = A^{-2}(x)$ takes the form
\[
V_0(x) \sim C_1^{-1} -c_2(x-1)^{2m_2 + 1}, \,\, m_2 \geq 1
\]
where $C_1>1$ and $c_2>0$.  Of course the constants are arbitrary
(chosen to agree with those in \cite{ChMe-lsm}), and $c_2$ could be
negative without changing much of the analysis.  This Lemma and the
proof are in \cite{ChMe-lsm}.

\begin{lemma}
\label{L:ml-inv-2}
For $\epsilon>0$ sufficiently small, let $\phi \in \s(T^* \reals)$
have compact support in $\{ |(x-1,\xi) |\leq \epsilon\}$.  Then there
exists $C_\epsilon>0$ such that 
\begin{equation}
\label{E:ml-inv-2}
\| P(z,h) \phi^w u \| \geq C_\epsilon
h^{(4m_2+2)/(2m_2+3)}  \|
\phi^w u \|, \,\,\, z \in [C_1^{-1}-\epsilon, C_1^{-1} + \epsilon].
\end{equation}
\end{lemma}

\begin{remark}
This estimate is also known to be sharp in the sense that the exponent
$(4m_2+2)/(2m_2+3)$ cannot be improved (see \cite{ChMe-lsm}).
\end{remark}

\subsubsection{Unstable infinitely degenerate and cylindrical critical elements}
\label{SS:unst-inf}

In this subsection, we study the case where the principal part of the
potential $V(x) = A^{-2}(x) + h^2 V_1(x)$ has
an infinitely degenerate maximum, say, at the point $x = 0$.  Let
$V_0(x) = A^{-2}(x)$.  As
usual, we again assume that $V_0(0) = 1$, so that
\[
V_0(x) = 1 - \O(x^\infty)
\]
in a neighbourhood of $x = 0$.  Of course this is not very precise, as
$V_0$ could be constant in a neighbourhood of $x = 0$ and still satisfy
this.
So let us first assume that $V_0(0) = 1$, and
$V_0'(x)$ vanishes to infinite order at $x = 0$, however, $\pm V_0'(x) <0$
for $\pm x >0$.  That is, the critical point at $x = 0$ is infinitely
degenerate but isolated.


\begin{lemma}
\label{L:ml-inv-3a}
For $\epsilon>0$ sufficiently small, let $\phi \in \s(T^* \reals)$
have compact support in $\{ |(x,\xi) |\leq \epsilon\}$.  Then for any
$\eta>0$, there
exists $C_{\epsilon,\eta}>0$ such that 
\begin{equation}
\label{E:ml-inv-3a}
\| P(z,h) \phi^w u \| \geq C_{\epsilon, \eta} {h^{2+\eta}} \|
\phi^w u \|, \,\,\, z \in [1-\epsilon, 1 + \epsilon].
\end{equation}
\end{lemma}

For our next result, we consider the case where there is a whole
interval at a local maximum value.  That is, we assume the principal part
of the effective potential
$V_0(x)$ has a maximum $V_0(x) \equiv 1$ on an interval, say $x \in
[-a,a]$, and that $\pm V_0'(x) <0$ for $\pm x > a$ in some neighbourhood.  

\begin{lemma}
\label{L:ml-inv-3b}
For $\epsilon>0$ sufficiently small, let $\phi \in \s(T^* \reals)$
have compact support in $\{ |x| \leq a + \epsilon,\, |\xi |\leq \epsilon\}$.  Then for any
$\eta>0$, there
exists $C_{\epsilon,\eta}>0$ such that 
\begin{equation}
\label{E:ml-inv-3b}
\| P(z,h) \phi^w u \| \geq C_{\epsilon, \eta} {h^{2+\eta}} \|
\phi^w u \|, \,\,\, z \in [1-\epsilon, 1 + \epsilon].
\end{equation}
\end{lemma}

\subsubsection{Infinitely degenerate and cylindrical inflection
  transmission critical elements}

\label{SS:inf-deg-infl}

In this subsection, we assume the effective potential has a critical
element 
 of infinitely degenerate or cylindrical inflection
transmission type.  This is very similar to Subsection \ref{SS:unst-inf},
but now the potential is assumed to be monotonic in a neighbourhood of
the critical value.

We begin with the case where the potential has an isolated
infinitely degenerate critical point of inflection transmission type.
As in the previous subsection, we write $V(x) = A^{-2}(x) + h^2
V_1(x)$ and denote $V_0(x) = A^{-2}(x)$ to be the principal part of
the potential.  
Let us assume the point $x = 1$ is
an infinitely degenerate inflection point, so that locally near $x = 1$,
the potential takes the form
\[
V_0(x) \sim C_1^{-1} - (x-1)^{\infty}, 
\]
where $C_1>1$.  Of course the constant is arbitrary
(chosen to again agree with those in \cite{ChMe-lsm}).  Let us assume
that our potential satisfies $V_0'(x) \leq 0$ near $x = 1$, with
$V_0'(x)<0$ for $x \neq 1$ in some neighbourhood so that the critical point $x = 1$ is isolated.

\begin{lemma}
\label{L:ml-inv-4a}
For $\epsilon>0$ sufficiently small, let $\phi \in \s(T^* \reals)$
have compact support in $\{ |(x-1,\xi) |\leq \epsilon\}$.  Then for
any $\eta>0$, there
exists $C = C_{\epsilon,\eta}>0$ such that 
\begin{equation}
\label{E:ml-inv-4a}
\| P(z,h) \phi^w u \| \geq C_\epsilon {h^{2 + \eta}} \|
\phi^w u \|, \,\,\, z \in  [C_1^{-1}-\epsilon, C_1^{-1} + \epsilon].
\end{equation}
\end{lemma}

On the other hand, if $V_0'(x) \equiv 0$ on an interval, say $x-1 \in
[-a,a]$ with $V_0'(x) < 0$ for $x-1 < -a$ and $x -1> a$, we do not expect
anything better than Lemma \ref{L:ml-inv-4a}.  The
next lemma says that this is exactly what we do get.  To fix an energy
level, assume $V_0 \equiv
C_1^{-1}$ on $[-a,a]$.  
\begin{lemma}
\label{L:ml-inv-4b}
For $\epsilon>0$ sufficiently small, let $\phi \in \s(T^* \reals)$
have compact support in $\{ |x-1| \leq a + \epsilon,/, |\xi| \leq \epsilon\}$.  Then for
any $\eta>0$, there
exists $C = C_{\epsilon,\eta}>0$ such that 
\begin{equation}
\label{E:ml-inv-4b}
\| P(z,h) \phi^w u \| \geq C_\epsilon {h^{2 + \eta}} \|
\phi^w u \|, \,\,\, z \in  [C_1^{-1}-\epsilon, C_1^{-1} + \epsilon].
\end{equation}
\end{lemma}

\section{Proof of Theorem \ref{T:T1} and Corollary \ref{C:C1}}

\label{S:proofs}



Recall the conjugated Laplacian is
\[
-\tDelta = - \p_x^2 - A^{-2}(x) \p_\theta^2 + V_1(x),
\]
where $V_1(x)$ has been computed above.  We will do some analysis and
reductions now {\it before} separating variables.  If we are
considering quasimodes
\[
(-\tDelta - \lambda^2) u = E(\lambda) \| u \|,
\]
where
\[
E(\lambda) = \O( \lambda^{-\beta_0})
\]
for some $\beta_0 >0$,
then we begin by rescaling.  Set $h = \lambda^{-1}$ so that
\[
(-h^2 \p_x^2 - h^2 A^{-2}(x) \p_\theta^2 + h^2 V_1(x) -1) u = \tE(h) \|
u \|,
\]
where $\tE(h) = h^2 E( h^{-1} ) = \O( h^{2 + \beta_0})$.  With $\xi,
\eta$ the dual variables to $x, \theta$ as usual, the semiclassical
symbol of this operator is
\[
p = \xi^2 + A^{-2}(x) \eta^2 + h^2 V_1(x)-1,
\]
and the semiclassical principal symbol is
\[
p_0 = \xi^2 + A^{-2}(x) \eta^2-1.
\]
It is worthwhile to point out that at this point our semiclassical
parameter is $h = \lambda^{-1}$.  After separating variables later in
the proof, we will let $h = |k|^{-1}$, where $k$ is the angular
momentum parameter.  However, in the regime where we so take $h$,
$|k|$ and $\lambda$ will be comparable, so it is merely a choice of convenience.

It is important to keep in mind for the remainder of this paper what the various parameters
represent.  Here, the variable $\eta$ represents $h D_\theta$.  As we will
eventually be decomposing in Fourier modes in the $\theta$ direction,
this means that the variable $\eta$ takes values in $h \ZZ$.

We next record that a standard $h$-parametrix argument tells us that
any quasimode is concentrated on the energy surface where $\{p_0 = 0
\}$.  The proof is standard.

\begin{lemma}
\label{L:char-surf}
Suppose $u$ satisfies
\[
(-h^2 \p_x^2 - h^2 A^{-2}(x) \p_\theta^2 + h^2 V_1(x) -1) u = \tE(h) \|
u \|,
\]
where $\tE(h) = h^2 E( h^{-1} ) = \O( h^{2 + \beta_0})$, and $\Gamma \in
\s^0$ satisfies $\Gamma \equiv 1$ in a small fixed neighbourhood of $\{p_0 = 0
\}$.  Then
\[
(1-\Gamma^w) u = \O( h^{2 + \beta_0} ).
\]
\end{lemma}

Hence we will restrict our attention to the characteristic surface
where $\{ p_0 = 0 \}$.
Using our moment map idea, we know that $\eta$ is invariant under the
classical flow.  Hence if $\eta$ is very large, our operator will be
elliptic, while if $\eta$ is very small, the parameter $\xi$ will be
bounded away from zero, and hence we will have uniform propagation
estimates.  Let us make this more precise.  Let $A_0 = \min (A(x))$
and $A_1 = \max (A(x))$, and let 
\[
1 = \psi_0(\eta) + \psi_1(\eta) + \psi_2(\eta)
\]
be a partition of unity satisfying
\[
\psi_0 \equiv 1 \text{ on } \{ | \eta |^2 \leq \frac{1}{2} A_0^2 \}
\]
with support in $\{ | \eta |^2 \leq \frac{3}{4} A_0^2 \}$; 
\[
\psi_2 \equiv 1 \text{ on } \{ | \eta |^2 \geq {2} A_1^2 \}
\]
with support in $\{ | \eta |^2 \geq \frac{3}{2} A_1^2 \}$.
Then, on $\supp \psi_0$, we have
\[
\eta^2 A^{-2}(x) \leq \eta^2 A_0^{-2} \leq \frac{3}{4},
\]
and on $\supp \psi_2$, we have
\[
\eta^2 A^{-2}(x) \geq \eta^2 A_1^{-2} \geq \frac{3}{2}.
\]

Now for our quasimode $u$, write
\[
u = u_0 + u_1 + u_2 + u_3 := \psi_0^w \Gamma^w u + \psi_1^w \Gamma^w u
+ \psi_2^w \Gamma^w u + (1 - \Gamma^w)u.
\]
Since $hD_\theta$ commutes with $-\tDelta$ and we can choose $\Gamma =
\Gamma(p_0)$ so that $[p_0^w, \Gamma^w] = \O(h^3)$, each of these $u_j$ are
also quasimodes of the same order as $u$ (but of course may have small
or even trivial $L^2$ mass).  

\subsubsection{Estimation of $u_0$}
\label{SSS:u0}
Observe that on the support of $\psi_0$, since $\eta$ is invariant, we
have $| \xi |^2 \geq 1/4 - \O(h^2)$, which means the propagation speed in the
$x$-direction is bounded below.  We claim this implies
\[
\| u_0 \|_{L^2_{x, \theta} } \leq c_0 \| u_0 \|_{L^2([a,b]_x \times
  \SS_\theta ) }
\]
for some $c_0 >0$.  In other words, $u_0$ is {\it uniformly}
distributed in the sense that the mass cannot be vanishing in $h$ on
any set.

The claim  follows by propagation of
singularities.  
The standard propagation of singularities result applies whenever the
classical flow propagates singularities from one region to another in
phase space.  Since we are analyzing the region where $\xi \neq 0$, we
have uniform propagation in the $x$ direction.  
A general statement is given in the following Lemma (a refinement of
H\"ormander's original result \cite{Hor-sing}).  For a proof in this context, see,
for example,  \cite[Lemma 6.1]{Chr-NC} and \cite[Lemma 4.1]{BuZw-bb}.

\begin{lemma}
\label{wf-lemma}
Suppose $V_0 \Subset T^*X$, $p$ is a symbol, $T >0$, $A$ an operator, and $V \Subset T^*X$ a neighbourhood of $\gamma$ satisfying
\begin{eqnarray}
\left\{ \begin{array}{l}
\forall \rho \in \{ p_0^{-1}(0) \} \setminus V, \,\,\, \exists \, 0 <t<T \,\, \text{and} \,\, \epsilon = \pm 1 \,\,\, 
\text{such that} \\
\exp(\epsilon s H_{p_0})(\rho) \subset \{ p_0^{-1}(0) \} \setminus V \,\, \text{for} \,\, 0 < s < t, \,\, \text{and} \\
\exp(\epsilon t H_{p_0})(\rho) \in V_0; \\
\end{array} \right.
\end{eqnarray}
and $A$ is microlocally elliptic in $V_0 \times V_0 $.  If $B \in \Psi^{0,0}(X, \Omega_X^\half)$ and $\WF (B) \subset T^*X 
\setminus V$, then
\begin{eqnarray}
\label{E:prop-00}
\left\| Bu \right\| \leq C \left( h^{-1} \left\| Pu \right\| + \| Au \| \right) + \O (h^\infty) \|u\|.
\end{eqnarray}
\end{lemma}

Fix two non-empty intervals in the $x$ direction, $(a,b)$ and
$(c,d)$ and assume $u = u_0$ is $L^2$ normalized.    
Now using that $Pu = \O(h^{2 + \beta_0}) \| u \|$, we have 
\begin{align*}
\| u \|_{L^2((c,d)\times \SS^1} & \leq C h^{-1} \|P u \| + C_2 \| u
\|_{L^2((a,b)\times \SS^1)} \\
& \leq C h^{1 + \beta_0} \| u \|_{L^2(\SS^1 \times \SS^1)} + C_2 \| u
\|_{L^2((a,b)\times \SS^1)},
\end{align*}
for some $C_2>0$.  
For $h>0$ sufficiently small, this implies if $u$ has mass bounded
below independent of $h$ in any $x$ neighbourhood $(c,d)$, then 
\[
\| u \|_{L^2((a,b)\times \SS^1)} \geq c'>0
\]
independent of $h$.  Rescaling in terms of $u_0$ if $u_0$ is not normalized, we recover
\[
\| u_0 \|_{L^2((a,b)\times \SS^1)} \geq c' \| u_0 \|.
\]

Since the interval $(a,b)$ is arbitrary, we have shown that the
$L^2$-mass on any rotationally invariant neighbourhood is positive independent of $h$.  Thus
\eqref{E:lower-bound} holds with a lower bound independent of $h = \lambda^{-1}$.


\subsubsection{Estimation of $u_2$}
On the other hand, on the support of $\psi_2$, we have the principal
symbol satisfies
\[
| p_0 | \geq \frac{1}{2},
\]
so we claim that an  elliptic argument shows 
\[
\| u_2 \|_{L^2} = \O(h^\infty) \| u_2 \|_{L^2}.
\]

That is, since $| p_0 | \geq \frac{1}{2}$ on support of $\psi_2$,
there is an $h$-parametrix for $P$ there: there exists $Q$ such that 
\[
Q P \psi_2^w = \psi_2^w + \O(h^\infty),
\]
and further $Q$ has bounded $L^2$ norm.  
Hence
\begin{align*}
\| u_2 \| & = \| QP u_2 \| + \O(h^\infty) \| u_2 \| \\
& \leq C \| P u_2 \| + \O(h^\infty) \| u_2 \| \\
& = \O(h^{2 + \beta_0} ) \| u_2 \|.
\end{align*}
This implies $u_2 = \O(h^\infty)$.





\subsubsection{Estimation of $u_1$}
\label{SSS:u1}


In order to consider the final part $u_1$, which is microsupported
where all the critical points are, we will employ one further
reduction.  Since $u_1$ is microsupported in a region where $| \eta |$
is bounded between two constants, say, $a_0 \leq | \eta | \leq a_1$,
and $\eta = h k$ for some integer $k$, 
a priori the number of angular momenta $k$ in the wavefront set of
$u_1$ is comparable to $h^{-1}$.  We can do better than that.  Using
the semiclassical calculus, we will next show that there exists $k_0
\in \ZZ$ such that for any $\epsilon>0$, we have
\[
u_1 = \sum_{|k - k_0| \leq h^{-\epsilon} } e^{i k \theta} \phi_k(x) +
\O(h^\infty) \| u_1 \|.
\]
That is, we claim that the Fourier decomposition of $u_1$ can actually
only have $\O(h^{-\epsilon})$ non-trivial modes.  To prove this claim,
fix $k_0 \in \ZZ$ and any $\epsilon>0$, and choose a $k_1 \in \ZZ$
satisfying
\[
| k_1 - k_0 | \geq h^{-\epsilon}.
\]
We will show that we can decompose $u_1$ into (at least) two pieces with disjoint
microsupport, one near $hk_0$ and one near $hk_1$.  Evidently, these
two pieces correspond to different angular momenta $\eta$, so have
wavefront sets associated to different level sets of the moment map.
Of course, level sets sufficiently close (in an $h$-dependent set) may
contribute to a single irreducible quasimode, but the point is to
quantify how far away from a single level set one needs to go before
leaving the microsupport of an irreducible quasimode.  

In order to make this rigorous, let $\eta_j = h k_j$ for $j = 0, 1$,
and choose $\chi(r) \in \Ci_c ( \reals) $ satisfying
\[
\chi(r) \equiv 1 \text{ for } | r | \leq 1 ,
\]
with support in $\{ | r | \leq 2 \}$.  For $j = 0, 1$, let
\[
\chi_j( \eta, h) = \chi \left( \frac{ \eta - \eta_j }{h^{1 -
      \epsilon/2}} \right).
\]
As semiclassical symbols, the $\chi_j$ are in a harmless $h^{1/2 -
  \epsilon/4}$ calculus, and moreover they only depend on $\eta$ (not on
$\theta$) and commute with $-\tDelta$.  
On the support of each of the $\chi_j$, we have
\[
\left| \frac{ \eta - \eta_j }{h^{1 - \epsilon/2} } \right| = \left|
  \frac{ h k - h k_j }{h^{1 - \epsilon/2} }\right| = \left| \frac{ k -
      k_j }{h^{- \epsilon/2}} \right| \leq 2.
\]
This implies
\[
| k - k_j | \leq 2 h^{-\epsilon/2},
\]
so as $h \to 0+$, $\chi_0$ and $\chi_1$ have disjoint supports.  This means the functions $\chi_1 u_1$ and
$\chi_2 u_1$ have disjoint $h$-wavefront sets, so they are almost
orthgonal:
\[
\lll \chi_1 u_1, \chi_2 u_1 \rrr = \O(h^\infty).
\]
Hence if each of these functions has nontrivial $L^2$ mass, then $u$
was not an irreducible quasimode.

Finally, we analyse the function $u_1$, but spread
over at most $\O(h^{-\epsilon})$ Fourier modes.  
Throughout the remainder of this section, let $\lambda$ be large and
fixed.  Let us consider a single Fourier mode confined 
to a single angular momentum $k$. 

The case of $u_1$ is the most interesting case, as the microsupport of
$u_1$ contains all
of the critical elements.  
Now recalling again the separated equation \eqref{E:separation} with the
potential
\[
V(x) = A^{-2}(x) + h^2 V_1(x),
\]
let $A_0$ and $A_1$ again be the min/max respectively of $A(x)$.  Our
spectral parameter now is $z = h^2 \lambda^2$.  
We are localized where 
\[
\frac{1}{2} A_0^2 \leq (\lambda^{-1} k)^2 \leq 2 A_1^2,
\]
or
\[
\frac{1}{2} A_1^{-2} \leq z \leq 2 A_0^{-2}.
\]
This of course implies that $\lambda$ and $k$ are
comparable.  Let $(a,b) \subset \SS^1$ be a non-empty interval.  We
need to show that if $u$ is a weak irreducible quasimode, 
\[
( (hD)^2 + V(x) -z) u = \O( h^{2 + \beta_0}) \| u \|,
\]
with $\| u \|= 1$, 
then either $\| u \|_{L^2(a,b)} = \O(h^\infty) =
\O(\lambda^{-\infty})$, or $\| u \|_{L^2(a,b)} \geq C_\epsilon h^{1 + \epsilon}$
for any $\epsilon>0$.  Let us assume that $u$ is nontrivial so that
$\|u \|_{L^2(a,b)} \geq c h^N$ for some $N$.



There are a number of subcases to consider here.  We observe that,
according to Lemma \ref{wf-lemma}, we
can always microlocalize further to a set close to the energy level of
interest.  That is, for $P(z,h) = (hD)^2 + V(x) -z$, if $P(z,h) u =
\O(h^{2 + \beta_0})$, then if $\psi(r) \in \Ci_c( \reals)$ satisfies
$\psi \equiv 1$ for $r$ near $0$, we have for any $\delta>0$
\[
\psi^w ( ( \xi^2 + V(x) - z) / \delta ) u = u + \O(h^{2 + \beta_0} ).
\] 
For the rest of this section, we write $\psi^w$ for this energy cutoff.



{\bf Case 1:}  

\begin{figure}
\hfill
\centerline{\input{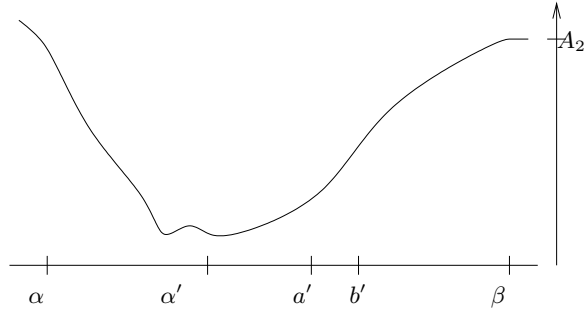}}
\caption{\label{fig:fig3b-1} The function $A^{-2}(x)$ and the weakly
  unstable critical point $\beta$.}
\hfill
\end{figure}

Next, assume $z$ is in a small neighbourhood of a critical energy level, and assume $A'(x)
\neq 0$ somewhere on $(a,b)$.  Then let $(a', b')$ be a non-empty
interval with 
\[
(a', b') \Subset \{ A' \neq 0\} \cap (a,b),
\]
and let $(\alpha', \beta) \supset (a', b')$ be the maximal connected
interval with $A'(x) \neq 0$ on $(\alpha', \beta)$.  Now $A'$ has constant sign on $(\alpha', \beta)$, so at least one of
$\alpha'$ or $\beta$ is part of a weakly unstable critical element
(see Figure \ref{fig:fig3b-1}).

Without loss in generality, assume $A' <0$ on $(\alpha', \beta)$ so
that at least $\beta$ lies in a weakly unstable critical element.
That is, the principal part of the potential $A^{-2}(x)$ increases as
$x \to \beta-$, and takes the value, say $A^{-2}(\beta) = A_2$.  Let
$(\alpha, \beta)$ be the maximal open interval containing $(\alpha',
\beta)$ where $A^{-2}(x) < A_2$ on $(\alpha, \beta)$.  As $A^{-2}(x) <
A^{-2}(\alpha)$ for $x \in (\alpha, \beta)$ and $A^{-2}(\alpha) = A_2$, we have $(A^{-2}(x) ) '
<0$ for $x \in (\alpha, \beta)$ sufficiently close to $\alpha$.  That
means that either $\alpha$ is part of a weakly unstable critical
element, or $A'(\alpha) \neq 0$.  We break the analysis into the two
separate subsubcases, beginning with $A'(\alpha) \neq 0$.

\begin{figure}
\hfill
\centerline{\input{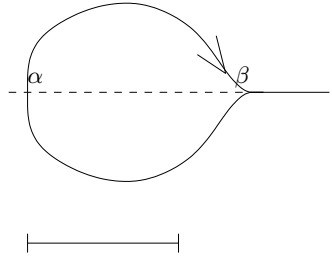}}
\caption{\label{fig:fig3b-2} If $A'(\alpha) \neq 0$, the unstable
  manifold from $\beta$ flows into the stable manifold at $\beta$
  (homoclinicity).  The interval indicates a region with propagation
  speed uniformly bounded below.}
\hfill
\end{figure}

{\bf Case 1a:}  
If $A'(\alpha) \neq 0$, then the weakly unstable/stable manifolds
associated to $(A^{-2})'(\beta) = 0$ 
are
homoclinic to each other (see Figure \ref{fig:fig3b-2}), and in particular, propagation of
singularities can be used to control the mass along this whole
trajectory, as long as we stay away from the right hand endpoint
$\beta$.  That is, propagation of singularities implies for any
$\eta>0$ independent of $h$, 
\begin{align*}
\| \psi^w u \|_{L^2(\alpha, \beta - \eta)} & \leq C_\eta (h^{-1} \|
((hD)^2 + V - z ) \psi^w u \| + \| \psi^w u \|_{L^2(a', b') } \\
& \leq C_\eta h^{1 + \beta_0} \| u \| + \| \psi^w u \|_{L^2(a', b') } .
\end{align*}
Hence by taking $h>0$ sufficiently small, we need to bound $\| \psi^w
u \|_{L^2(\alpha, \beta - \eta)}$ from below in terms of $\| u \|$.

Let
$[\beta, \kappa]$ be the maximal connected interval containing $\beta$
on which $A' = 0$ (we allow $\kappa = \beta$ if the critical point is isolated).  
Let $\tchi \equiv 1$ on $[\beta, \kappa]$ with
support in a small neighbourhood thereof, and let $\chi \equiv 1$ on
$\supp \tchi$ with support in a slightly smaller set so that $(1 -
\tchi) \geq (1 - \chi)$ and $(1 - \tchi) \geq c | \chi'|$.  Then
writing $P(z,h) = (hD)^2 + V - z$, we have from Theorem
\ref{T:inf-deg-est} (for any $\epsilon>0$)
\begin{align*}
\| u \| & \leq \| \chi u \| + \| ( 1 - \chi) u \| \\
& \leq C_\epsilon  h^{-2-\epsilon} \| P(z,h) \chi u \| + \| (1 - \tchi
) u \| \\
& \leq C_\epsilon h^{-2-\epsilon} (\| \chi P(z,h) u \| + \| [P(z,h) ,
\chi ] u \| )  + \| (1 - \tchi
) u \| \\
& \leq C_\epsilon' ( h^{\beta_0 - \epsilon} \| u \| +
h^{-1-\epsilon}\| (1 - \tchi) u \| )  + \| (1 - \tchi
) u \| 
\end{align*}
Rearranging and taking $h>0$ sufficiently small and $\epsilon <
\beta_0$, we get
\begin{equation}
\label{E:tchi-1}
\| (1 - \tchi) u \| \geq C_\epsilon h^{1 + \epsilon} \| u \|.
\end{equation}
Now either the wavefront set of $u$ is contained in the closure of the
lift of
$(\alpha, \beta)$ or it isn't.  In the latter case there is nothing to
prove.  In the former case, we conclude that $u = \O(h^\infty)$ on any
open subset whose closure does not meet the set $[\alpha, \beta]$.  We
appeal to propagation of singularities one more time.  Since
$A'(\alpha) \neq 0$, propagation of singularities applies in a
neighbourhood of $\alpha$, so that (shrinking $\eta>0$ if necessary)
for some $c_1>0$, 
\[
\| u \|_{L^2(\alpha, \beta - \eta)} \geq c_1 \| u \|_{L^2( \alpha -
  \eta, \beta - \eta)}.
\]
Since we have assumed $u = \O(h^\infty)$ on $(\alpha - \eta, \beta +
\eta)^c$, 
this estimate, together with \eqref{E:tchi-1} and \eqref{E:prop-00} allows us to conclude
\[
\| u \|_{L^2(\alpha, \beta - \eta)} \geq C \| (1 - \tchi) u \| \geq C_\epsilon h^{1 + \epsilon}
\| u \|.
\]

\begin{figure}
\hfill
\centerline{\input{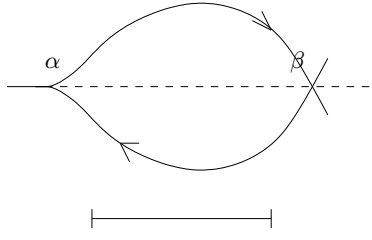}}
\caption{\label{fig:fig3b-3} If $A'(\alpha) = 0$, the unstable
  manifold from $\beta$ flows into the stable manifold at $\alpha$ and
  vice versa.  The interval indicates a region with propagation
  speed uniformly bounded below.}
\hfill
\end{figure}

{\bf Case 1b:}  
We now  consider
the possibility that $A'(\alpha) = 0$ as well as $A'(\beta) = 0$ (see Figure \ref{fig:fig3b-3}).  In this case,
propagation of singularities fails at both endpoints of $(\alpha,
\beta)$, so we can only conclude that for  any $\eta>0$
independent of $h$,
\[
\| u \|_{L^2(\alpha + \eta, \beta - \eta)} \leq C_\eta (h^{-1} \|
((hD)^2 + V - z ) u \| + \| u \|_{L^2(a', b') }.
\]

Hence now it suffices to prove that for some $\eta>0$ small but
independent of $h$, we have the estimate
\[
\| u \|_{L^2(\alpha + \eta, \beta - \eta)}  \geq C_\epsilon h^{1 +
  \epsilon} \| u \|
\]
for any $\epsilon>0$.

Let
$[\beta, \kappa]$ be the maximal connected interval containing $\beta$
on which $A' = 0$, and let $[\omega, \alpha]$ be the maximal connected
interval containing $\alpha$ on which $A' = 0$.  
Let $\tchi \equiv 1$ on $[\beta, \kappa] \cup [\omega, \alpha]$ with
support in  small neighbourhoods thereof, and let $\chi \equiv 1$ on
$\supp \tchi$ with support in a slightly smaller set so that $(1 -
\tchi) \geq (1 - \chi)$ and $(1 - \tchi) \geq c | \chi'|$.  Since both
$[\omega, \alpha]$ and $[\beta, \kappa]$ are weakly unstable, we can apply
Theorem \ref{T:inf-deg-est} and the same argument as above to finish
this case.

{\bf Case 2:}  
Finally, we assume $(a,b) \subset \{ A' = 0 \}$.
Again, if $A^{-2} \equiv A_3$ on $(a,b)$ and $z \neq A_3$, we can use
propagation of singularities to control $\| u \|_{L^2(a,b)}$ from
below by its mass on the connected component in $\{ p = z \}$
containing $(a,b)$ (as in the case of $u_0$ above).  Hence we are interested in the
case where $z$ is in a small neighbourhood of $A_3$.

If $u = \O(h^\infty)$ on $(a,b)$ there is nothing to prove, so assume
not.  Then if $[\alpha, \beta] \supset (a,b)$ is the
maximal connected interval where $A^{-2}(x) \equiv A_3$, the wavefront
set of $u$ is contained in a small neighbourhood of $[\alpha, \beta]$, so that for 
$\delta>0$ as small as we like by taking a sufficiently localized
energy cutoff, we have
\[
\| u \|_{L^2( [\alpha - \delta, \beta + \delta]^c )} = \O_\delta (h^\infty).
\]

\begin{figure}
\hfill
\centerline{\input{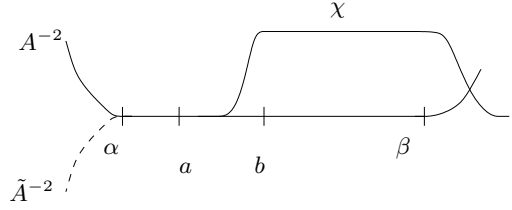}}
\caption{\label{fig:fig3c} The setup for Case 2.  Here if the
  quasimode is small in $(a,b)$, we cut off to the right of $(a,b)$
  and modify $A^{-2}$ to the left to be weakly unstable.  We then arrive at a contradiction.}
\hfill
\end{figure}

That means that, either
\[
\| u \|_{L^2([a,b])} \geq c >0, \,\,\, 
\| u \|_{L^2([\alpha - \delta, a])} \geq c >0, \text{ or }
\| u \|_{L^2([b, \beta +
  \delta])} \geq c >0.
\]
If the first estimate is true, we're done, so assume without loss in
generality that  $\| u \|_{L^2([b, \beta +
  \delta])} \geq c >0$.  Assume for contradiction that there exists
  $\epsilon_0>0$ such that $\| u \|_{L^2(a,b)} \leq C h^{1 +
    \epsilon_0}$.  Let $\chi \in \Ci_c$ be a smooth function such that
  $\chi \equiv 1$ on $[b, \beta + \delta]$ with support in $(a, \beta
  + 2 \delta)$.  Write $\tu = \chi u$.  If $[\alpha, \beta]$ is a
  weakly stable critical element, modify $A^{-2}(x)$ on the support of
  $1-\chi$ so that $[\alpha, \beta]$ is weakly unstable.  That is, if
  $(A^{-2}(x))' <0$ for $x < \alpha$ in some neighbourhood, replace
  $A$ with a locally defined function $\tA$ satisfying $\tA \equiv A$
  on $\supp \chi$ but $(\tA^{-2}(x))' >0$ for $x < \alpha$ in some
  neighbourhood.  If $[\alpha, \beta]$ is weakly unstable, then let
  $\tA \equiv A$ (see Figure \ref{fig:fig3c}.  We apply Theorem \ref{T:inf-deg-est} once again
  (for any $\epsilon>0$):
\begin{align*}
\| \chi u \| & \leq C_\epsilon h^{-2 - \epsilon} \| ((hD)^2 + \tA^{-2} + h^2 V_1 -z )
\chi u \| \\
& = C_\epsilon h^{-2 - \epsilon} \| ((hD)^2 + A^{-2} + h^2 V_1 -z )
\chi u \| \\
& \leq C_\epsilon h^{-2 - \epsilon} ( \| P(z,h) u \| + \| [P(z,h),
\chi] u \| ) \\
& \leq C_\epsilon'  (h^{\epsilon_0 - \epsilon}  \| u \| + h^{-1 -
  \epsilon} \| u \|_{L^2(a,b)} ) + \O(h^\infty),
\end{align*}
where the $\O(h^\infty)$ error comes from the part of the commutator
$[P(z,h), \chi]$ supported outside a neighbourhood of $[\alpha,
\beta]$ (the other part contributing the integral over $(a,b)$).  But
our contradiction assumption implies that the right hand side is
$o(1)$ as $h \to 0$ provided $\epsilon< \epsilon_0$.  As $\| \chi u \|
\geq c >0$, this is a contradiction.

\subsection{Finishing up the proof}

We now put together the estimates of $u_0, u_1, u_2, u_3$. 
Since $u_3 = \O( h^{2 + \beta_0} )$ and $u_2 = \O(h^\infty)$, for
$h>0$ sufficiently small,  at least one of $u_0$ and $u_1$ must have $L^2$
mass bounded below independent of $h$.  If $u_0$ has $L^2$ mass
bounded below independent of $h$ we're done by the propagation of
singularities argument in Subsection \ref{SSS:u0}.  Hence we need to
conclude Theorem \ref{T:T1} assuming $u_0$ is small and $u_1$ carries most of
the $L^2$ mass.

Fix $(a,b)$ as considered in Subsection \ref{SSS:u1} and recall we
know that for any $\epsilon>0$ 
\[
u_1 = \sum_{|k - k_0| \leq h^{-\epsilon} } e^{i k \theta} \phi_k(x) +
\O(h^\infty) \| u_1 \|.
\]
We use the notation $\Omega = (a,b)_x \times \SS^1_\theta$ as in the
statement of Theorem \ref{T:T1}.  
Each $\phi_k$ satisfies either
\[
\| \phi_k \|_{L^2(a,b)} = \O(|k|^{-\infty})
\]
or for any $\epsilon>0$, 
\[
\| \phi_k \|_{L^2 (a,b) } \geq c_2 |k|^{-1-\epsilon} \| \phi_k \|_{L^2
  (\SS^1_x)}.
\]
In the first case, these $\phi_k$s have disjoint wavefront sets from
the $\phi_k$s in the latter case, so leaving them in the sum would
mean our   quasimode was not irreducible.  Removing these from the sum
and reindexing if necessary, we conclude
\begin{align*}
\| u_1 \|_{L^2(\Omega)}^2 & = \sum_{|k - k_0| \leq k_0^{\epsilon} }  \|\phi_k(x)\|_{L^2(a,b)}^2 +
\O(h^\infty) \| u_1 \|^2  \\
& \geq c_2'  k_0^{-2-2\epsilon} \sum_{|k - k_0| \leq k_0^{\epsilon} }
\|\phi_k(x)\|_{L^2(\SS^1_x )}^2 -
\O(h^\infty) \| u_1 \|^2 \\
& = c_2'' \lambda^{-2-2\epsilon} \| u_1 \|_{L^2(\SS^1_x \times \SS^1_\theta)}^2 - \O(
\lambda^{-\infty})  \| u_1 \|^2.
\end{align*}
This concludes the proof of Theorem \ref{T:T1}.

\qed

\bibliographystyle{alpha}
\bibliography{ef-surf-rev-bib}

\end{document}